\documentclass[11pt, reqno, a4paper]{article}
\usepackage{amsfonts}
\usepackage{mathrsfs}
\usepackage{amssymb}
\usepackage{amsmath}
\newtheorem{theorem}{Theorem}[section]
\newtheorem{lemma}{Lemma}[section]
\newtheorem{corollary}{Corollary}[section]
\newtheorem{definition}{Definition}[section]
\newtheorem{remark}{Remark}[section]

\begin{document}

\title{The Direct and Converse Inequalities for Jackson-Type Operators
on Spherical Cap
\thanks{The research was supported by the National
Natural Science Foundation of China (No. 60873206) and the Natural
Science Foundation of Zhejiang Province of China (No. Y7080235)}}

\author{ Yuguang Wang \and
 Feilong Cao\thanks{Corresponding
author: Feilong Cao,  E-mail: \tt flcao@263.net;
feilongcao@gmail.com}}

  \date{}
\maketitle

\begin{center}
\footnotesize

 Department of Information and Mathematics Sciences,
 China Jiliang University,

 Hangzhou, 310018, Zhejiang Province,  P. R. China.

\begin{abstract}
Approximation on the spherical cap is different from that on the
sphere which requires us to construct new operators. This paper
discusses the approximation on the spherical cap. That is, so called
Jackson-type operator $\{J_{k,s}^m\}_{k=1}^{\infty}$ is constructed to
approximate the function defined on the spherical cap $D(x_0,\gamma)$.
We thus establish the direct and inverse inequalities and obtain
saturation theorems for $\{J_{k,s}^m\}_{k=1}^{\infty}$ on the cap $D(x_0,\gamma)$.
 Using methods of $K$-functional and multiplier, we
obtain the inequality
\begin{eqnarray*}
C_1\:\| J_{k,s}^m(f)-f\|_{D,p}\leq
\omega^2\left(f,\:k^{-1}\right)_{D,p} \leq C_2 \max_{v\geq k}\|
J_{v,s}^m(f) - f\|_{D,p}
\end{eqnarray*}
and that the saturation order of these operators is $O(k^{-2})$,
where $\omega^2\left(f,\:t\right)_{D,p}$ is the modulus of smoothness of degree 2,
the constants $C_1$ and $C_2$ are independent of $k$ and $f$.

{\bf MSC(2000):}  41A17

{\bf Key words:} approximation; saturation order; Jackson-type
operator; spherical cap
\end{abstract}
\end{center}

\section{Introduction}
\markright{Course Paper of China Jiliang University}
In the past decades, many mathematicians dedicated to establish the
Jackson and Bernstein-type theorems on the sphere. Early works, such as Butzer
\cite{Butzer1971}, Nikol'ski\v{i} \cite{Nikolskii1984}
and
\cite{Nikolskii1987}, Lizorkin \cite{Lizorkin1983} had
successfully established the direct and inverse theorems on the
sphere. In early 1990s, Li and Yang \cite{Li1991} constructed Jackson operators
on the sphere and obtained the Jackson and Bernstein-type theorems for the Jackson operators.

 Jackson operator on the sphere
 is defined by (see \cite{Lizorkin1983}, \cite{Li1991})
\begin{eqnarray}\label{eq23}
J_{k,s}(f)(x):= \frac{1}{|\mathbb{S}^{n-2}|} \int_{\mathbb{S}^{n-1}}
f(y)\mathscr{D}_{k,s}(\arccos x\cdot y)f(y)d\omega(y)
\end{eqnarray}
where $k$ and $s$ are positive integers,
$$
\mathscr{D}_{k,s}(\theta)
=\left({ \frac{\sin\frac{k\theta}{2}}
{\sin\frac{\theta}{2}}}\right)^{2s}
$$
 is the
classical Jackson kernel, $f$ is measurable function of degree $p$
on the sphere $\mathbb{S}^{n-1}$ in $\mathbb R^n$, $d \omega(y)$ is
the elementary surface piece, $|\mathbb{S}^{n-1}|$ is the
measurement of $\mathbb{S}^{n-1}$. For $f\in L^p(\mathbb{S}^{n-1})$,
$(1\leq p\leq\infty)$ ($L^\infty(\mathbb{S}^{n-1})$ is the
collection of continuous functions on $\mathbb{S}^{n-1}$), Li and
Yang \cite{Li1991} proved that
\[C_1\:\|J_{k,s}(f)-f\|_{{\mathbb S^{n-1}},p}\leq\omega^2(f,k^{-1})\leq
C_2\:k^{-2}\sum_{v=1}^k v\|J_{v,s}(f)-f\|_{{\mathbb S^{n-1}},p}\]
and the saturation
order for $J_{k,s}$ is $k^{-2}$, where $C_1$ and $C_2$ are
independent of positive integer $k$ and $f$, and
$\omega^2(f,t)$
is
the modulus of smoothness of degree 2 on the unit sphere
${\mathbb S^{n-1}}$.

On the
spherical caps, we desire to obtain the same results. The main difficulty we face is
to establish the Bernstein-type inequality (for polynomials) on the cap.
Recently, Belinsky, Dai and Ditzian \cite{Belinsky2003} constructed
$m$-th translation operator $S^{m}_{\theta}$ when discussing the averages of functions on the sphere.
This inspires us and allows us to construct the $m$-th
Jackson-type operator $J_{k,\;s}^{m}$ on the spherical cap. Fortunately, the Bernstein-type
inequality  holds for $J_{k,\;s}^{m}$,
which helps us get the direct and inverse inequalities of approximation on the spherical cap.
Finally, we obtain that the saturation order for the constructed Jackson-type
operator is $k^{-2}$, the same to that of the Jackson operator on the sphere.

\section{Definitions and Auxiliary Notations}

Throughout out this paper, we  denote by the letters $C$ and
$C_i$ ($i$ is either positive integers or variables on which $C$
depends only) positive constants depending only on the dimension
$n$. Their value may be different at different occurrences, even
within the same formula.
 We shall denote the points in
 $\mathbb {S}^{n-1}$
 by $x, x_0,x_1,\dots$,
 $y, y_0,y_1,\dots$,
 and
 $z,z_0,z_1,\dots$,
and the elementary surface piece on $\mathbb {S}^{n-1}$ by
$d\omega$. If it is necessary, we shall write $d\omega(x)$ referring
to the variable of the integration. The notation $a \approx b$ means
that there exists a positive constant $C$ such that $C^{-1}b\leq a
\leq Cb$ where $C$ is independent of some variable $n$ on which $a$
and $b$ both depend.

Next, we introduce some concepts and properties of sphere as well as
caps (see \cite{Muller1966}, \cite{Wang2000}).
 The volume of
$\mathbb{S}^{n-1}$ is
$$
\Omega_{n-1}:=\int_{\mathbb{S}^{n-1}}d\omega= \frac{2\pi^{n/2}}{\Gamma(n/2)}.
$$
Corresponding to $d\omega$,  the inner product on
$\mathbb{S}^{n-1}$ is defined by
$$
\left\langle f,g\right\rangle:=\int_{\mathbb{S}^{n-1}}f(x)\overline{g(x)}d\omega(x).
$$
Denote by
$L^p(\mathbb{S}^{n-1})$
 the space of
$p$\:-integrable functions on $\mathbb{S}^{n-1}$ endowed with the
norms
$$
\| f\|_{\infty}:=\|f\|
_{L^\infty(\mathbb{S}^{n-1})}:=\mbox{ess}\!\!\sup_{x\in
\mathbb{S}^{n-1}}|f(x)|
$$
and
$$
\|f\|_p:=\|f\|_{L^p(\mathbb{S}^{n-1})}:=\left\{\int_{\mathbb{S}^{n-1}}|f(x)|^pd\omega(x)\right\}^{1/p}<\infty,
\quad 1\leq p<\infty.
$$

We denote by $D(x_0,\gamma)$ the spherical cap with center  $x_0$
and angle $0<\gamma\leq\frac{\pi}{2}$, i.e.,
$$
D(x_0,\gamma):=\left\{x\in \mathbb{S}^{n-1}:x\cdot
x_0\geq\cos\gamma\right\},
$$
and by $D(\gamma)$ the volume of
$D(x_0,\gamma)$, i.e.,
$$D(\gamma):=\int_0^\gamma|\mathbb{S}^{n-2}|\sin^{n-2}\theta d\theta.
$$
Then for fixed $x_{0}$ and $\gamma$, $L^{p}(D(x_{0},\gamma))$ is a
Banach space endowed with the norm $\|\cdot\|_{D,p}$
defined by
$$
\| f\|_{D,\infty}:=\mbox{ess}\!\!\!\!\!\sup_{x\in
D(x_0,\gamma)}\!\!\!\!\!|f(x)|
$$
and
$$
\|f\|_{D,p}\::=\:\left\{\int_{D(x_0,\gamma)}|f(x)|^pd\omega(x)\right\}^{1/p},\;1\leq
p<\infty.
$$
 For any $f\in
L^{p}(D(x_{0},\gamma))$, we note
$$
f^*(x)=\left\{\begin{array}{ll}
f(x), & x\in D(x_0,\gamma),\\
0, & x\in \mathbb{S}^{n-1}\backslash D(x_0,\gamma),
\end{array}
\right.
$$
and clearly, $f^*\in L^{p}(\mathbb{S}^{n-1})$ and
$\|f^*\|_{p}=\|f\|_{D,p}$. This allows us to introduce some operators
on spherical cap using existing operators on the sphere.

\begin{definition}
Suppose $T: L^{p}(\mathbb{S}^{n-1})\rightarrow
L^{p}(\mathbb{S}^{n-1})$ is an operator on $\mathbb{S}^{n-1}$, then
$$T_{x_{0},\gamma}: L^{p}(D(x_0,\gamma))\rightarrow L^{p}(D(x_0,\gamma))$$
$$T_{x_{0},\gamma}(f)(x)= T(f^{*})(x),\quad x\in D(x_0,\gamma)$$  is
called the operator on $D(x_0,\gamma)$  introduced by T. We may
use the notation $T$ instead of $T_{x_{0},\gamma}$ for
convenience  if without mixing up.
\end{definition}

 We now make a brief introduction
of projection operators $Y_{j}(\cdot)$ by ultraspherical
(Gegenbauer) polynomials $\{G_{j}^{\lambda}\}_{j=1}^{\infty}$
 $\left(\lambda={\frac{ n-2}{2}}\right)$ for discussion
of saturation property of Jackson operators.

Ultraspherical polynomials $\{G_j^{\lambda}\}_{
j=1}^{\infty}$ are defined in terms of the generating function (see
\cite{Stein1971}):
\[\frac{1}{(1-2 t r+r^2)^\lambda}=\sum_{j=0}^{\infty}G_j^{\lambda}(t)r^j\]
where $|r|<1$, $|t|\leq 1$.

For any $\lambda>0$, we have (see \cite{Stein1971})
\begin{equation}\label{eq24}
G_1^{\lambda}(t)=2\lambda t
\end{equation} and
\begin{equation}\label{eq27}
 \frac{d}{dt}G_j^{\lambda}(t)= 2\lambda
G_{j-1}^{\lambda+1}(t).
\end{equation}
When $ \lambda=\frac{n-2}{2}$ (see \cite{Wang2000}),
\[G_j^{\lambda}(t)=\frac{\Gamma(2\lambda+j)}{\Gamma(j+1)\Gamma(2\lambda)}P_j^n(t),\quad
j=0,1,2,\dots\] where $P_j^n(t)$ is the Legendre polynomial of
degree $j$.
Particularly,
\[G_j^{\lambda}(1)=\frac{\Gamma(2\lambda+j)}{\Gamma(j+1)\Gamma(2\lambda)}P_j^n(1)=\frac{\Gamma(2\lambda+j)}{\Gamma(j+1)\Gamma(2\lambda)},\quad
j=0,1,2,\dots.\]
Therefore,
\begin{equation}\label{eq28}
P_j^n(t)=\frac{G_j^{\lambda}(t)}{G_j^{\lambda}(1)}.
\end{equation}
Besides, for any $j=0,1,2,\ldots$, and $|t|\leq 1$,
$|P_j^n(t)|\leq1$ (see \cite{Muller1966}).

The projection operators is defined by
\begin{eqnarray*}
Y_j(f)(x)=\frac{\Gamma(\lambda)(n+\lambda)}{2\pi^{\frac{n}{2}}}\int_{\mathbb{S}^{n-1}}G_j^{\lambda}(x\cdot
y)f^{*}(y)\:d\omega(y).
\end{eqnarray*}
It follows from (\ref{eq24}), (\ref{eq27}) and (\ref{eq28}) that
\begin{eqnarray}\label{eq25}
\lim_{t\rightarrow 1-}
{ \frac{1-(P_j^n(t))^m}{1-(P_1^n(t))^m}=\frac{j(j+2\lambda)}{2\lambda+1}}.
\end{eqnarray}

In the same way, we define the inner product  on
$D(x_0,\gamma)$ as follows,
$$
\left\langle
f,g\right\rangle_D:=\int_{D(x_0,\gamma)}f(x)\overline{g(x)}d\omega(x).
$$

We denote by $\widetilde{\Delta}$ the Laplace-Beltrami operator
$$
\widetilde{\Delta} f:=\sum_{i=1}^n\frac{\partial^2g(x)}{\partial
x_i^2}\bigg|_{|x|=1},\quad g(x)=f\left(\frac{x}{|x|}\right),
$$
by which we  define a $K$-function on $D(x_0,\gamma)$ as
$$
K(f,\delta)_{D,p}:=\inf\left\{\|f-g\|_{D,p}+\delta\|\widetilde{\Delta}
g\|_{D,p}:\ g,\; \widetilde{\Delta} g\in L^p(D(x_0,\gamma))\right\}.
$$
For $f\in L^1(D(x_0,\gamma))$,  the translation operator is defined
by
$$
S_\theta(f)(x):=\frac1{|\mathbb{S}^{n-2}|\sin^{n-2}\theta}\int_{x\cdot
y=\cos \theta}f^*(y)d\omega'(y),
$$
where $d\omega'(y)$ denotes the the elementary surface piece on the
sphere $\{y\in D(x_0,\gamma):x\cdot y=\cos \theta\}$. Then we have
\begin{eqnarray*}
\int_{D(x_0,\gamma)}f(x)d\omega(x)=\int_0^\pi
S_\theta(f)(x_0)|\mathbb{S}^{n-2}|\sin^{n-2}\theta d\theta.
\end{eqnarray*}
The modulus of smoothness of $f$  is defined by
$$
\omega^2(f,\delta)_{D,p}:=\sup_{0<\theta\leq \delta}\|S_\theta (f)-f\|_{D,p}.
$$
Using the method of \cite{Butzer1971}, we have
\begin{eqnarray}\label{eq22}
C_1\omega^2(f,\delta)_{D,p}\leq K(f,\delta^2)_{D,p}\leq
C_2\omega^2(f,\delta)_{D,p}.
\end{eqnarray}

We introduce $m$-th translation operator in terms of multipliers
(see \cite{Lizorkin1987}, \cite{Rustamov1994}, \cite{Wang2000})
\begin{eqnarray}\label{eq271}
S^{m}_{\theta}(f)= \sum^{\infty}_{j=0}\left(\frac{G^{\lambda}_{j}(\cos\theta)}{G^{\lambda}_{j}(1)}\right)^{m}Y_{j}(f)
\equiv\sum^{\infty}_{j=0}\left(P^n_{j}(\cos\theta)\right)^{m}Y_{j}(f),\;
f\in L^p(D(x_0,\gamma)).
\end{eqnarray}
It has been proved that (see \cite{Wang2000})
\begin{eqnarray}\label{eq26}
S_\theta(f)(x) =
\frac1{|\mathbb{S}^{n-2}|\sin^{n-2}\theta}\int_{x\cdot y
=\cos\theta}f^{*}(y)d\omega'(y) & = &
\sum^{\infty}_{j=0}P^n_{j}(\cos\theta)Y_{j}(f)(x)\nonumber\\
& = &S^{1}_{\theta}(f)(x)
\end{eqnarray}

With the help of $S_{\theta}^{m}$, we can construct Jackson-type
operators on $D(x_0,\gamma)$.

\begin{definition}\label{def22}
For $f\in L^p(D(x_0,\gamma))$, the $m$-th Jackson-type operator of
degree $k$ on $D(x_0,\gamma)$ is defined by
$$J_{k,s}^{m}(f)(x)=\int_{0}^{\gamma}S_{\theta}^{m}(f)(x)\widetilde{D}_{k,s}(\theta)\sin^{2\lambda}\theta d\theta,$$
where $ \lambda=\frac{n-2}{2}$, and
$\widetilde{D}_{k,s}(\theta)=A_{\gamma,k,s}^{-1}\left({ \frac{\sin^{2s}\frac{k\theta}{2}}{\sin^{2s-1}\frac{\theta}{2}}}\right)$
satisfying
$ \displaystyle\int_{0}^{\gamma}\widetilde{D}_{k,s}(\theta)\sin^{2\lambda}\theta
d\theta=1.$
\end{definition}

\begin{remark}
We may notice that $\widetilde{D}_{k,s}(\theta)$ is a bit different
from classical Jackson kernel
$$\mathscr{D}_{k,s}(\theta)=\left({ \frac{\sin\frac{k\theta}{2}}{\sin\frac{\theta}{2}}}\right)^{2s}.$$
This difference will help us to prove Bernstein inequality for
$J_{k,s}^{m}$. For sake of ensuring that Bernstein inequality for
$J_{k,s}^{m}$ holds, $\gamma$ has to be no more than
$ \frac{\pi}{2}$. Particularly, for $m=1$, we have
\begin{eqnarray*}
J_{k,s}^{1}(f)(x) & = &
\int_{0}^{\gamma}S_{\theta}^{1}(f)(x)\widetilde{D}_{k,s}(\theta)\sin^{2\lambda}\theta
d\theta\\
& = &
\int_0^\gamma\frac{1}{|\mathbb{S}^{n-2}|\sin^{2\lambda}\theta}\int_{x\cdot
y
=\cos\theta}f^{*}(y)d\omega'(y)\widetilde{D}_{k,s}(\theta)\sin^{2\lambda}\theta
d\theta\\
& = &
\frac{1}{|\mathbb{S}^{n-2}|}\int_{D(x_0,\gamma)}f^{*}(y)\widetilde{D}_{k,s}(\arctan(x\cdot
y))d\omega(y)\\
& = & \frac{1}{|\mathbb{S}^{n-2}|}\int_{D(x_0,\gamma)\bigcap
D(x,\gamma)}f(y)\widetilde{D}_{k,s}(\arctan(x\cdot y))d\omega(y).
\end{eqnarray*}
\end{remark}

Finally, we introduce the definition of saturation for operators
(see \cite{Berens1968}).
\begin{definition}
 Let $\varphi(\rho)$ be a positive function with respect to $\rho$,
$0<\rho<\infty$, tending monotonely to zero as
$\rho\rightarrow\infty$. For $\rho>0$, $I_\rho$ is a sequence of
operators. If there exists $\mathcal {K}\subseteqq
L^p(D(x_0,\gamma))$ such that:\vspace{0.1 cm}\\\parbox{6
cm}{\begin{tabular}{ll}
(i) & If\; $\| I_\rho(f)-f\|_{D,p} = o(\varphi(\rho))$, then $I_\rho(f)=f$; \\
(ii) & $\| I_\rho(f)-f\|_{D,p} = O(\varphi(\rho))$ if and only if
$f\in\mathcal {K}$;
\end{tabular}}\parbox{1 cm}{}\vspace{0.2 cm}\\
then $I_\rho$ is said to be saturated on $L^p(D(x_0,\gamma))$ with
order $O(\varphi(\rho))$ and $\mathcal {K}$ is called its saturation
class.
\end{definition}

\section{Some Lemmas}

In this section, we show some lemmas on both $S_{\theta}^{m}$ and
$J_{k,s}^{m}$ as the preparation for the main results. For
$S_{\theta}^{m}$, we have
\begin{lemma}\label{lm31}
For $f\in L^{p}(D(x_{0},\gamma))$, $1\leq p\leq\infty$,
$0<\theta\leq\pi$,\\\parbox{6 cm}{\begin{tabular}{ll} (i) & For
$1\leq m_{1}< m$,
$S_{\theta}^{m}=S_{\theta}^{m_{1}}S_{\theta}^{m-m_{1}}$;\\
 (ii)  & For\; $m\geq 1$, $\| S_{\theta}^{m}(f)\|_{D,p}\ \leq \| f\|_{D,p}$;\\
 (iii) & For\; $m\geq 1$, $\| S_{\theta}^{m}(f)-f\|_{D,p}\leq m \| S_{\theta}(f)-f\|_{D,p}$;\vspace{0.1 cm}\\
 (iv)  &  For\; $m > \frac{2([\frac{n}{2}]+3)}{n-2}$, $0<\theta\leq \frac{\pi}{2}$, $\|
\widetilde{\Delta}S_{\theta}^{m}(f)\|_{D,p}\leq C_m\theta^{-2}\| f\|_{D,p}$,\\
& where $C_m\rightarrow 0$, as $m\rightarrow \infty$;\\
 (v)   & For\; $m\geq 1$, and $f$ which satisfies $\widetilde{\Delta}f \in
L^{p}(D(x_0,\gamma))$,
$\|\widetilde{\Delta}S_{\theta}^{m}(f)\|_{D,p}\leq \|
\widetilde{\Delta}f \|_{D,p}$.
\end{tabular}}\parbox{1 cm}{}\\
\end{lemma}
{\bf Proof.} $(i)$, $(ii)$ and $(iii)$ are clear. Using Remark 3.5
of \cite{Belinsky2003}, we can obtain $(iv)$. For $(v)$, we have
\begin{eqnarray*} \widetilde{\Delta}S_{\theta}(f)=
-\sum_{j=0}^{\infty}j(j+2\lambda)P_j^n(\cos\theta)Y_{j}(f)=
S_{\theta}(\widetilde{\Delta}f),
\end{eqnarray*}
which implies
\begin{eqnarray*}
\|\widetilde{\Delta}S_{\theta}^{m}(f)\|_{D,p} \:=\: \|
S_{\theta}\left(\widetilde{\Delta}S_{\theta}^{m-1}(f)\right)\|_{D,p}
\:\leq \:  \| \widetilde{\Delta}S_{\theta}^{m-1}(f)\|_{D,p} \leq
\cdots
 \:\leq\:  \| \widetilde{\Delta}f\|_{D,p}.\quad  \Box
\end{eqnarray*}

We need the following lemma.
\begin{lemma}\label{lm32}
For $\beta\geq-2$, $2s>\beta+2\lambda+1$,
$ 0<\gamma\leq \pi$, $n\geq3$, and $k\geq1$, we
have
\begin{eqnarray}\label{eq34}
\int_0^\gamma \theta^\beta
\widetilde{D}_{k,s}(\theta)\sin^{2\lambda}\theta\:d\theta\approx
k^{-\beta},
\end{eqnarray}
where ${ \lambda=\frac{n-2}{2}}$.
\end{lemma}
{\bfseries Proof.} A simple calculation gives, for $\beta\geq-2$,
\begin{eqnarray*}
\int_{0}^{\gamma}\theta^{\beta}\frac{\sin^{2s}\frac{k\theta}{2}}{\sin^{2s-1}\frac{\theta}{2}}\sin^{2\lambda}\theta
d\theta \leq C \left(\int_{0}^{\infty}\frac{\sin^{2s}}{\theta^{2s-(2\lambda+\beta+1)}} d\theta\right) k^{2s-(2\lambda+\beta+2)},
\end{eqnarray*}
\begin{eqnarray*}
\int_{0}^{\gamma}\theta^{\beta}\frac{\sin^{2s}\frac{k\theta}{2}}{\sin^{2s-1}\frac{\theta}{2}}\sin^{2\lambda}\theta
d\theta \geq C \left(\int_{0}^{\frac{k\gamma}{4}}\frac{\sin^{2s}}{\theta^{2s-(2\lambda+\beta+1)}} d\theta\right) k^{2s-(2\lambda+\beta+2)},
\end{eqnarray*}
\begin{eqnarray*}
\int_{0}^{\gamma}\theta^{\beta}\widetilde{D}_{k,s}(\theta)\sin^{2\lambda}\theta
d\theta =
\frac{ \displaystyle\int_{0}^{\gamma}\theta^{\beta}\frac{\sin^{2s}\frac{k\theta}{2}}{\sin^{2s-1}\frac{\theta}{2}}\sin^{2\lambda}\!\theta
d\theta}{ \displaystyle\int_{0}^{\gamma}\frac{\sin^{2s}\frac{k\theta}{2}}{\sin^{2s-1}\frac{\theta}{2}}\sin^{2\lambda}\!\theta d\theta}
 \approx k^{-\beta}.\quad\Box
\end{eqnarray*}

For Jackson-type operators, we have
\begin{lemma}\label{lm33}
For $f\in L^{p}(D(x_{0},\gamma))$, $1\leq p<\infty$,\vspace{0.2 cm}\\
\parbox{6 cm}{\begin{tabular}{ll}
(i)  & $\|J_{k,s}^{m}(f)\|_{D,p}\leq\| f\|_{D,p}$;  \\
(ii) & For $f$ which satisfies $\widetilde{\Delta}f \in
L^{p}(D(x_0,\gamma))$, $\|
\widetilde{\Delta}J_{k,s}^{m}(f)\|_{D,p}\leq \|
\widetilde{\Delta}f\|_{D,p}$;\vspace{0.15 cm}\\
(iii)& For $n \geq 3$, $m
>  \frac{2\left([\frac{n}{2}]+3\right)}{n-2}$ and
$0<\gamma\leq \frac{\pi}{2}$,
$\|\widetilde{\Delta}J_{k,s}^{m}(f)\|_{D,p}\leq C_mk^{2}\|
f\|_{D,p}$.
\end{tabular}}\parbox{1 cm}{}
\end{lemma}
{\bf Proof.} From the definition and $(ii)$ and $(v)$ of
Lemma~\ref{lm31}, $(i)$ and $(ii)$ are clear. We just have to add
the proof of $(iii)$. In fact, using Minkowski inequality, $(iv)$ of
Lemma~\ref{lm31} and Lemma~\ref{lm32}, we have
\begin{eqnarray*}
\| \widetilde{\Delta}J_{k,s}^{m}(f)\|_{D,p} &\leq   &
\int_{0}^{\gamma}\|
\widetilde{\Delta}S_{\theta}^{m}(f)\|_{D,p}\widetilde{D}_{k,s}(\theta)\sin^{2}\theta
d\theta\\
 &\leq   &  C_m \|
f\|_{D,p}\int_{0}^{\gamma}\theta^{-2}\widetilde{D}_{k,s}(\theta)\sin^{2\lambda}\theta
d\theta \\
&\approx&  C_{m,\gamma}\; k^{2}\| f\|_{D,p},
\end{eqnarray*}
where the constant in the approximation is independent of $m$ and
$k$.\quad$\Box$

The following lemma is useful in the proof of Bernstein-type
inequality for Jackson-type operators.

\begin{lemma}
[\cite{Wickeren1986}]
\label{lm39} Suppose that for nonnegative sequences
$\left\{\sigma_k\right\}_{k=1}^{\infty}$,
$\left\{\tau_k\right\}_{k=1}^{\infty}$ with $\sigma_1=0$ the
inequality
\begin{eqnarray*}
\sigma_n\leq\left(\frac{k}{n}\right)^p\sigma_k+\tau_k,\quad
p>0,\quad 1\leq k\leq n
\end{eqnarray*}
is satisfied for any positive integer $n$. Then one has
\begin{eqnarray}
\sigma_n\leq C_p n^{-p}\sum_{k=1}^{n} k^{p-1}\tau_k.
\end{eqnarray}
\end{lemma}

The following lemma gives the multiplier representation of
$J_{k,s}^m(f)$, which follows from Definition~\ref{def22} and
(\ref{eq271}).
\begin{lemma}\label{lm37}
For $f\in L^{p}(D(x_0,\gamma))$, $J_{k,s}^m(f)$ has the
representation
\begin{eqnarray}
J_{k,s}^m(f)(x) & = & \sum_{j=0}^{\infty}\xi_k^m(j)Y_{j}(f)(x)
\end{eqnarray}
where \[\xi_k^m(j)=\int_0^{\gamma}
\widetilde{D}_{k,s}(\theta)\left(P_j^n(\cos\theta)
\right)^m\sin^{2\lambda}\theta\:d\theta,\quad j=0,1,2,\dots.\]
\end{lemma}

The following lemma is useful for determining the saturation order.
It can be deduced by the methods of \cite{Berens1968} and
\cite{Butzer1972}.
\begin{lemma}\label{lm38}
Suppose that $\left\{I_\rho\right\}_{\rho>0}$ is a sequence of
operators on $L^p(D(x_0,\gamma))$, and there exists function series
$\left\{\lambda_{\rho}(j)\right\}_{j=1}^{\infty}$ with respect to
$\rho$, such that
\[I_\rho(f)(x)=\sum_{j=0}^{\infty}\lambda_\rho(j)Y_j(f)(x)\]for every $f\in L^p(D(x_0,\gamma))$. If for any $j=0,1,2,\dots$,
 there $\mathit{exists}$ $\varphi(\rho)\rightarrow
0\!+(\rho\rightarrow\rho_0)$ such that
\[\lim_{\rho\rightarrow\rho_0}\frac{1-\lambda_\rho(j)}{\varphi(\rho)}=\tau_j\neq0,\] then
$\left\{I_\rho\right\}_{\rho>0}$ is saturated on
$L^p(D(x_0,\gamma))$ with the order $O(\varphi(\rho))$ and the
collection of all constants is the invariant class for
$\left\{I_\rho\right\}_{\rho>0}$ on $L^p(D(x_0,\gamma))$.
\end{lemma}

\section{Main Results and Their Proofs}

In this section, we shall discuss the main results, that is, the
lower and upper bounds as well as the saturation order for
Jackson-type operators on $L^p(D(x_0,\gamma))$.

The following theorem gives the Jackson-type inequality for
$J_{k,s}^m$.

\begin{theorem}\label{th1}
For any integer $m\geq 1$ and $0<\gamma\leq \frac{\pi}{2}$, $\{J_{k,s}^m\}_{k=1}^{\infty}$
is the series of Jackson-type operators on $L^p(D(x_0,\gamma))$ defined above, and $g\in L_p^{(2)}(D(x_0,\gamma)):=\{f\in
L^p(D(x_0,\gamma)):\widetilde{\Delta} f\in
L^p(D(x_0,\gamma))\},\;1\leq p\leq \infty$. \\Then
\begin{eqnarray*}
\|J_{k,s}^m(g)-g\|_{D,p}\leq C_{m,\gamma}\:k^{-2}\|\widetilde{\Delta}
g\|_{D,p}.
\end{eqnarray*}
Therefore, for $f\in L^p(D(x_0,\gamma))$,
\begin{eqnarray}\label{eq42}
\|J_{k,s}^m(f)-f\|_{D,p}\leq C_{m,\gamma}\:
\omega^2\left(f,k^{-1}\right)_{D,p},
\end{eqnarray}
where $C$ is independent of $k$ and $f$.
\end{theorem}
{\bfseries Proof.} Since $g\in L_p^{(2)}(D(x_0,\gamma)),$
  we have (see \cite{Berens1968})
\begin{eqnarray}\label{eq410}
S_\theta(g)(x)-g(x)=\int_0^\theta\sin^{2-n}\nu\int_0^\nu\sin^{n-2}\tau
S_\tau(\widetilde{\Delta} g)(x)d\tau d\nu.
\end{eqnarray}
and it is true that
\begin{eqnarray*}
\sup_{\theta>0}\left\{\theta^{-2}\int_0^\theta\sin^{2-n}\nu\int_0^\nu\sin^{n-2}\tau
d\tau d\nu\right\}<\infty.
\end{eqnarray*}
Therefore (explained below),
\begin{eqnarray*}
&&\|J_{k,s}^m(g)-g\|_{D,p}
 =
 \left\|\int_0^{\gamma}
\widetilde{D}_{k,s}(\theta)(S_\theta^m(g)(x)-g(x))\sin^{n-2}\theta d\theta\right\|_{D,p}\\
&\leq&
 m\:\int_0^\gamma
 \widetilde{D}_{k,s}(\theta)\sin^{n-2}\theta\left(\int_0^\theta\sin^{2-n}\nu\int_0^\nu\sin^{n-2}\tau
\|S_\tau(\widetilde{\Delta} g)\|_{D,p}d\tau d\nu\right)d\theta\\
&\leq&
m\:\sup_{\theta>0}\left\{\theta^{-2}\int_0^\theta\sin^{2-n}\nu\int_0^\nu\sin^{n-2}\tau
d\tau d\nu\right\} \|\widetilde{\Delta}
g\|_{D,p}\int_0^\gamma\theta^2\widetilde{D}_{k,s}(\theta)\sin^{n-2}\theta\:
d\theta\\
&\leq& C_{m,\gamma}\:k^{-2}\:\|\widetilde{\Delta} g\|_{D,p},
\end{eqnarray*}
where the Minkowski inequality, (\ref{eq410}) and Lemma~\ref{lm31}
are used in the first inequality, and the second one and third one
is deduced from (\ref{eq34}) and Lemma~\ref{lm31}. From (\ref{eq22})
and $(i)$ of Lemma~\ref{lm33}, it is easy to deduce
(\ref{eq42}).\quad$\Box$

Next, we prove the Bernstein-type inequality for $J_{k,s}^m(f)(x)$
for $f\in L^p(D(x_0,\gamma))$.

\begin{theorem}\label{th2}
Assume $\{J_{k,s}^m\}_{k=1}^{\infty}$, $m>
\frac{2\left([\frac{n}{2}]+3\right)}{n-2}$ are $m$-th Jackson-type
operators on
 $D(x_0,\gamma)$. For $f\in L^p(D(x_0,\gamma))$, $0<\gamma\leq \frac{\pi}{2}$, then there exits a constant $C$ independent of
$k$ and $f$ such that
\begin{eqnarray}\label{eq415}
\omega^2(f, \;k^{-1})_{D,p}\leq C \max_{v\geq k} \| J_{v,s}^m(f) - f
\|_{D,p}
\end{eqnarray}
holds for every $f\in L^p(D(x_0,\gamma))$ and every integer k.
\end{theorem}

\textbf{Proof.} Li and Yang \cite{Li1991} have proved the
Marchaud-Ste\v{c}kin inequality for Jackson operator on the sphere.
Following the method in \cite{Li1991}, we first prove the
Marchaud-Ste\v{c}kin inequality for $J_{k,s}^m$:
\begin{eqnarray}\label{eq416}
\omega^2(f, \:k^{-1})_{D,p}\leq C_{1} k^{-2} \sum_{v=1}^k v \|
J_{v,s}^m(f) - f \|_{D,p}.
\end{eqnarray}
Let
\[\sigma_{k}=k^{-2}\|\widetilde{\Delta} J_{k,s}^m(f)\|_{D,p}, \quad
\tau_k=\| J_{k,s}^m(f)-f \|_{D,p},\]then for $v=1,2,\dots,k$, by
Lemma~\ref{lm33},
\begin{eqnarray*}
\sigma_{k} & = & k^{-2}\|\widetilde{\Delta}
J_{k,s}^m(f)-\widetilde{\Delta}
J_{k,s}^m(J_{v,s}^m(f))+\widetilde{\Delta} J_{k,s}^m(J_{v,s}^m(f))\|_{D,p}\\
           & = & k^{-2}\|\widetilde{\Delta} J_{k,s}^m(f-J_{v,s}^m(f))\|_{D,p}+\|\widetilde{\Delta}
           J_{k,s}^m(J_{v,s}^m(f))\|_{D,p}\\
           & \leq & k^{-2}\left(C\:k^2\|
           f-J_{v,s}^m(f)\|_{D,p}+\|\widetilde{\Delta}
           J_{v,s}^m(f)\|_{D,p}\right)\\
           & = & C\:\|
           f-J_{v,s}^m(f)\|_{D,p}+k^{-2}\|\widetilde{\Delta}
           J_{v,s}^m(f)\|_{D,p}\\
           & = & C\:\tau_{v}+\left(\frac{v}{k}\right)^2\sigma_{v},
\end{eqnarray*}
so we can deduce from Lemma~\ref{lm39} (where $p$ is set to be $2$)
that\[\sigma_{k}\leq C\:k^{-2}\sum_{v=1}^{k}v\tau_v\] that is,
\begin{eqnarray}
\|\widetilde{\Delta} J_{k,s}^m(f)\|_{D,p}\leq C\:\sum_{v=1}^{k}v\|
J_{v,s}^m(f)-f\|_{D,p}.
\end{eqnarray}
 Since there exists ${ \frac{k}{2}}\leq r\leq k$ such that
\[\| J_{r,s}^m(f)-f \|_{D,p}=\min_{\frac{k}{2}\leq
v\leq k}\| J_{v,s}^m(f)-f\|_{D,p}\] then,
\begin{eqnarray*}
 K(f,\: k^{-2})_{D,p}
 &\leq& \| J_{r,s}^m(f)-f\|_{D,p}+k^{-2} \|\widetilde{\Delta}
J_{r,s}^m(f)\|_{D,p}\\
 &\leq& 4k^{-2} \sum_{\frac{k}{2}\leq v\leq k }v\|
 J_{v,s}^m(f)-f\|_{D,p}+C k^{-2}\sum_{v=1}^{k}v\|
 J_{v,s}^m(f)-f\|_{D,p}\\
 &\leq& C \:k^{-2} \sum_{v=1}^{k}v\|
 J_{v,s}^m(f)-f\|_{D,p}.
\end{eqnarray*}
By (\ref{eq22}), we obtain that
\begin{eqnarray*}
\omega^2(f,k^{-1})_{D,p}
 &\leq& C K(f,k^{-2})_{D,p}\\
 &\leq& C k^{-2} \sum_{v=1}^{k}v\|
 J_{v,s}^m(f)-f\|_{D,p}.
\end{eqnarray*}
So (\ref{eq416}) holds, and it implies that
\begin{eqnarray}\label{eq419}
\omega^2(f, k^{-1})_{D,p}\leq C_{1} k^{-1-\frac{1}{2}} \sum_{v=1}^k
v^{\frac{1}{2}} \| J_{v,s}^m(f) - f \|_{D,p}.
\end{eqnarray}
\indent In order to prove (\ref{eq415}), we have to show that
\begin{eqnarray}\label{eq420}
 \omega^2(f, k^{-1})_{D,p} &\approx& \frac{1}{k^2} \max_{1 \leq v \leq k}
 v^2 \| J_{v,s}^m(f) - f \|_{D,p}\nonumber\\
 &\approx& \frac{1}{k^{2+\frac{1}{4}}} \max_{1 \leq v \leq k}
 v^{2+\frac{1}{4}} \| J_{v,s}^m(f) - f \|_{D,p}.
\end{eqnarray}
We first prove
\begin{eqnarray}\label{eq421}
\omega^2(f, k^{-1})_{D,p} \approx \frac{1}{k^2} \max_{1 \leq v \leq k} v^2
\| J_{v,s}^m(f) - f \|_{D,p}.
\end{eqnarray}
It follows from (\ref{eq42}) and (\ref{eq419}) that
\begin{eqnarray*}
 \omega^2(f, k^{-1})_{D,p}
 &\leq& C_{1}k^{-2}k^{\frac{1}{2}}\sum_{v=1}^k
 v^{-\frac{3}{2}} \bigg(v^2 \| J_{v,s}^m(f) - f \|_{D,p}\bigg) \\
 &\leq& C_{1}k^{-2}\bigg(\max_{1\leq v \leq
 k} v^2 \| J_{v,s}^m(f) - f \|_{D,p}\bigg) \bigg(k^{\frac{1}{2}} \sum_{v=1}^k v^{-\frac{3}{2}}\bigg)\\
 &\leq& C_{3}k^{-2}\max_{1\leq v \leq k} v^2 \| J_{v,s}^m(f) - f \|_{D,p} \\
 &\leq& C_{4}k^{-2}\max_{1\leq v \leq k} v^2 \:\omega^2(f, v^{-1})_{D,p} \\
 &\leq& C_{4}\bigg(k^{-2}\max_{1\leq v \leq k} v^2\bigg(1+\left(\frac{k}{v}\right)^2\bigg)\bigg)\omega^2(f, \:k^{-1})_{D,p} \\
 &\leq& 2C_{4}\omega^2(f, k^{-1})_{D,p}.
\end{eqnarray*}
Then we prove
\begin{eqnarray*}
\omega^2(f, k^{-1})_{D,p} \approx \frac{1}{k^{2+\frac{1}{4}}} \max_{1
\leq v \leq k} v^{2+\frac{1}{4}} \| J_{v,s}^m(f) - f \|_{D,p}.
\end{eqnarray*}
In fact, the proof is similar to that of (\ref{eq421}).
\begin{eqnarray*}
\omega^2(f,k^{-1})_{D,p}
 &\leq& C_{1}k^{-2-\frac{1}{4}}k^{\frac{3}{4}}\sum_{k=1}^k
 v^{-\frac{7}{4}} \bigg(v^{2+\frac{1}{4}} \| J_{v,s}^m(f) - f \|_{D,p}\bigg) \\
 &\leq& C_{1}k^{-2-\frac{1}{4}}\bigg(\max_{1\leq v \leq k}
 v^{2+\frac{1}{4}} \| J_{v,s}^m(f) - f \|_{D,p}\bigg) \bigg(k^{\frac{3}{4}} \sum_{v=1}^k v^{-\frac{7}{4}}\bigg)\\
 &\leq& C_{5}k^{-2-\frac{1}{4}} \max_{1\leq v \leq k}
 v^{2+\frac{1}{4}} \| J_{v,s}^m(f) - f \|_{D,p} \\
 &\leq& C_{6}k^{-2-\frac{1}{4}} \max_{1\leq v \leq k} v^{2+\frac{1}{4}} \omega^2(f, v^{-1})_{D,p} \\
 &\leq& C_{6}\bigg(k^{-2-\frac{1}{4}}\max_{1\leq v \leq k} v^{2+\frac{1}{4}}\bigg(1+\left(\frac{k}{v}\right)^2\bigg)\bigg)\omega^2(f, k^{-1})_{D,p} \\
 &\leq& 2C_{6}\omega^2(f, k^{-1})_{D,p}.
\end{eqnarray*}
Hence,
\begin{eqnarray*}
 \omega^2(f, k^{-1})_{D,p} \approx \frac{1}{k^2} \max_{1 \leq v \leq k} v^2
 \| J_{v,s}^m(f) - f \|_{D,p}
 \approx \frac{1}{k^{2+\frac{1}{4}}} \max_{1 \leq v \leq k}
 v^{2+\frac{1}{4}} \| J_{v,s}^m(f) - f \|_{D,p}.
\end{eqnarray*}

Now we can complete the proof of (\ref{eq415}). Let
\begin{eqnarray*}
\max_{1 \leq v \leq k}
 v^{2+\frac{1}{4}} \| J_{v,s}^m(f) - f \|_{D,p}=
 k_{1}^{2+\frac{1}{4}} \| J_{k_{1},s}^m(f) - f
 \|_{D,p},\quad 1\leq k_{1}\leq k.
\end{eqnarray*}
it follows from (\ref{eq420}) that
\begin{eqnarray*}
 k^{-2} k_{1}^2\| J_{k_{1},s}^m(f) - f \|_{D,p}
 &\leq& k^{-2}\max_{1\leq v \leq k} v^2 \| J_{v,s}^m(f) - f \|_{D,p}\\
 &\leq& \frac{C_{7}}{k^{2+\frac{1}{4}}}\max_{1\leq v \leq k} v^{2+\frac{1}{4}} \| J_{v,s}^m(f) - f \|_{D,p}\\
 &=& \frac{C_{7}}{k^{2+\frac{1}{4}}} k_{1}^{2+\frac{1}{4}} \| J_{k_{1},s}^m(f) - f \|_{D,p}.
\end{eqnarray*}

Thus, $C_{7}^{-4} k \leq k_{1}$. Since $k_{1} \leq k$, then $k
\approx k_{1}$, we obtain from (\ref{eq420}) that
\begin{eqnarray*}
\omega^2(f,k^{-1})_{D,p}
 &\leq& C_{8}k^{-2-\frac{1}{4}}\max_{1 \leq v \leq k} v^{2+\frac{1}{4}} \| J_{v,s}^m(f) - f \|_{D,p} \\
 &=& C_{8}k^{-2-\frac{1}{4}} \bigg(k_{1}^{2+\frac{1}{4}} \| J_{k_{1},s}^m(f) - f \|_{D,p}\bigg)\\
 &\leq& C_{8}k^{-2-\frac{1}{4}} \max_{k_{1} \leq v \leq k} v^{2+\frac{1}{4}} \| J_{v,s}^m(f) - f \|_{D,p} \\
 &\leq& C_{8}\max_{k_{1} \leq v \leq k} \| J_{v,s}^m(f) - f \|_{D,p}.
\end{eqnarray*}
Noticing that $k \approx k_{1}$, we may rewrite the above inequality
as
\begin{eqnarray*}
\omega^2(f,\: k^{-1})_{D,p}\leq C_{9} \max_{v\geq k} \| J_{v,s}^m(f) -
f \|_{D,p}.
\end{eqnarray*}
This completes the proof.\quad$\Box$

We thus obtain the corollary of Theorem~\ref{th1} and Theorem~\ref{th2}
that:

\begin{corollary}
Suppose that $\{J_{k,s}^m\}_{k=1}^{\infty}$, $m>
 \frac{2\left([\frac{n}{2}]+3\right)}{n-2}$, are
Jackson-type operators on the spherical cap $D(x_0,\gamma)$,
$0<\gamma\leq\frac{\pi}{2}$, then the following are
equivalent for any $f\in$
$L^{p}\left(D(x_{0},\:\gamma)\right)$, $\alpha >0$,\\
\parbox{3 cm}{\begin{tabular}{ll} (i)&$\| J_{k,s}^m(f) -
f\|_{D,p}=O(k^{-\alpha}),$\quad$k\rightarrow \infty$;  \\
(ii)&$\omega^2(f,\delta)_{D,p}=O(\delta^{\:\alpha}),\quad\delta\rightarrow 0+$.\quad $\Box$
\end{tabular}} \hfill \parbox{1 cm}{}
\end{corollary}

\begin{theorem}
Suppose that $\{J_{k,s}^m\}_{k=1}^{\infty}$, $m\geq1$ are
Jackson-type operators on the spherical cap $D(x_0,\gamma)$,
$ 0<\gamma\leq\frac{\pi}{2}$. Then
$\{J_{k,s}^m\}_{k=1}^{\infty}$ are saturated on
$L^{p}\left(D(x_{0},\:\gamma)\right)$ with order $k^{-2}$ and the
collection of constants is their invariant class.
\end{theorem}

\textbf{Proof.}
We obtain from Lemma~\ref{lm32} that, for $v=1,2,\dots$,
\begin{eqnarray*}
 1-\xi_k^m(1)
  =
 \int_0^{\gamma}\widetilde{D}_{k,s}(\theta)\left(1-\left(\frac{G_1^\lambda(\cos\theta)}{G_1^\lambda(1)}\right)^m\right)\sin^{2\lambda}\theta\:d\theta
 & \approx &
 \int_0^{\gamma}\widetilde{D}_{k,s}(\theta)\sin^2\frac{\theta}{2}\sin^{2\lambda}\theta\:d\theta\\
 & \approx & k^{-2}.
\end{eqnarray*}
By Lemma~\ref{lm38}, if it is true that for
$j=0,1,2,\dots$,
\begin{eqnarray*}
\lim_{k\rightarrow\infty}\frac{1-\xi_k^m(j)}{1-\xi_k^m(1)} & = &
\frac{j(j+2\lambda)}{2\lambda+1},
\end{eqnarray*}
then the proof is completed.

In fact, for any $0<\delta<\gamma$, it follows from (\ref{eq34})
that
\begin{eqnarray*}
\int_{\delta}^{\gamma}\widetilde{D}_{k,s}(\theta)\sin^{2\lambda}\theta\:d\theta
& \leq &
\int_{\delta}^{\gamma}{\left(\frac{\theta}{\delta}\right)}^{\!\!3}\widetilde{D}_{k,s}(\theta)\sin^{2\lambda}\theta\:d\theta\\
& \leq &
\delta^{-3}\int_{0}^{\gamma}\theta^{\:3}\widetilde{D}_{k,s}(\theta)\sin^{2\lambda}\theta\:d\theta
\:\leq\; C_{\delta,s}k^{-3}.
\end{eqnarray*}

We deduce from (\ref{eq25}) that, for any $\varepsilon>0$, there exists
$\delta>0$, for $0<\theta<\delta$, it holds that
\begin{eqnarray*}
\left|\left(1-\left(P_j^n(\cos\theta)\right)^m\right)-\frac{j(j+2\lambda)}{2\lambda+1}\left(1-\left(P_1^n(\cos\theta)\right)^m\right)\right|
& \leq & \varepsilon\left(1-\left(P_1^n(\cos\theta)\right)^m\right).
\end{eqnarray*}
Then it follows that
\begin{eqnarray*}
& & \left|\left(1-\xi_k^m(j)\right)-\frac{j(j+2\lambda)}{2\lambda+1}\left(1-\xi_k^m(1)\right)\right|\\
& = &
\left|\int_{0}^{\gamma}\widetilde{D}_{k,s}(\theta)\left(1-\left(P_j^n(\cos\theta)\right)^m\right)\sin^{2\lambda}\theta\:d\theta\right.\\
& &
\left.-\int_{0}^{\gamma}\widetilde{D}_{k,s}(\theta)\left(1-\left(P_1^n(\cos\theta)\right)^m\right)\frac{j(j+2\lambda)}{2\lambda+1}\sin^{2\lambda}\theta\:d\theta\right|\\
& = &
\left|\int_{0}^{\gamma}\widetilde{D}_{k,s}(\theta)\left(\big(1-\left(P_j^n(\cos\theta)\right)^m\big)-
\big(1-\left(P_1^n(\cos\theta)\right)^m\big)\frac{j(j+2\lambda)}{2\lambda+1}\right)\sin^{2\lambda}\theta\:d\theta\right|\\
& \leq &
\int_{0}^{\delta}\widetilde{D}_{k,s}(\theta)\:\varepsilon\:\sin^{2\lambda}\theta\:d\theta
+2 \int_{\delta}^{\gamma}\widetilde{D}_{k,s}(\theta)\sin^{2\lambda}\theta\left(1+\frac{j(j+2\lambda)}{2\lambda+1}\right)\:d\theta\\
& \leq & \varepsilon k^{-2}+C_{\delta,s}k^{-3}.
\end{eqnarray*}
 So,
\begin{eqnarray*}
{ \lim_{k\rightarrow\infty}\frac{1-\xi_k(j)}{1-\xi_k(1)}=\frac{j(j+2\lambda)}{2\lambda}\neq
0}.
\end{eqnarray*}

Therefore, we obtain  by Lemma~\ref{lm38} that the saturation order
for $J_{k,s}^m$ is $O(1-\xi_k(1))\approx k^{-2}$.\quad $\Box$

\end{document}